\documentstyle[12pt]{article}

\input{mssymb}

\newtheorem{theorem}{\sc Theorem}[section] 
\newtheorem{lemma}[theorem]{\sc Lemma} 
\newtheorem{claim}[theorem]{\sc Claim} 
\newtheorem{subclaim}[theorem]{\sc Subclaim} 
\newtheorem{proposition}[theorem]{\sc Proposition}
\newtheorem{Observation}[theorem]{\sc Observation}

\newcommand{\gk}{\kappa}
\newcommand{\gl}{\lambda}
\newcommand{\ga}{\alpha}
\newcommand{\gTh}{\Theta}
\newcommand{\gf}{\phi}
\newcommand{\gr}{\rho}
\newcommand{\go}{\omega}
\newcommand{\gth}{\theta}
\newcommand{\gd}{\delta}
\newcommand{\gz}{\zeta}
\newcommand{\gb}{\beta}
\newcommand{\grg}{\gamma}
\newcommand{\gm}{\mu}
\newcommand{\gn}{\nu}
\newcommand{\gs}{\sigma}
\newcommand{\ha}{\aleph}
\newcommand{\dmd}{\diamondsuit}
\newcommand{\set}[2]{\{#1:#2\}}
\newcommand{\cf}{{\rm cf\/}}
\newcommand{\und}{*}
\newcommand{\Brch}{{\rm Brch}}
\newcommand{\se}{\subseteq}
\newcommand{\rest}{\mathord{\restriction}}
\newcommand{\proof}{{\sc Proof}\ } 
\newcommand{\fin}{\hspace{0.2in}\vrule width 6pt height 6pt depth 0pt 
\vspace{0.1in}}

\title{Some Compact Logics  - Results in ZFC }
\author{Alan H. Mekler\thanks{Research supported by NSERC}, Simon
Fraser University \and Saharon Shelah\thanks{Publication
\#375. Research supported by the BSF}, The Hebrew University\\ Rutgers
University\\
{\bf Dedicated to the memory of Alan by his friend, Saharon}}

\newcommand{\gint}[1]{[{\gamma^-_{#1,i}}, {\gamma^+_{#1,i}})}
\newcommand{\mon}[1]{{\cal #1}} 
\newcommand{\slt}{\prec^\otimes_\gth}
\newcommand{\dns}{does not split over}
\newcommand{\her}[1]{{\rm H}(#1)}
\newcommand{\Ob}{\mbox{\rm Ob}(\gd)}
\newcommand{\Qb}{Q_{\rm Brch}}
\newcommand{\Qba}{Q_{\rm Ba}}
\newcommand{\Qof}{Q_{\rm Of}}
\newcommand{\swd}{somewhere dense}
\newcommand{\nwd}{nowhere dense}
\newcommand{\vare}{\varepsilon}
\renewcommand{\iff}{if and only if}
\font\eufb=eufb10
\newcommand{\GC}{\hbox{\eufb\char"43}}

\begin{document}
\maketitle

\newpage

\section{Preliminaries}

While first order logic has many nice properties it lacks
expressive power. On the other hand 
second order logic is so strong
that it fails to have nice model theoretic properties such as
compactness. It is desirable to find natural logics which are stronger
than the first order logic but which still satisfy the compactness
theorem. Particularly attractive are those logics which allow
quantification over natural algebraic objects. One of the most natural
choices is to quantify over automorphisms of a structure (or
isomomorphisms between substructures). Generally compactness fails badly
\cite{Sh56}, but if we restrict ourself to certain concrete classes then
we may be able to retain compactness. In this paper we will show that if
we enrich first order logic by allowing quantification over isomorphisms
between definable ordered fields the resulting logic, $L(\Qof)$, is fully
compact. In this logic, we can give standard compactness proofs of
various results. For example, to prove that there exists arbitrarily
large rigid real closed fields, fix a cardinal $\gk$ and form the
$L(\Qof)$ theory in the language of ordered fields together with $\gk$
constants which says that the constants are pairwise distinct and the
field is a real closed field which is rigid. (To say the field is rigid
we use the expressive power of $L(\Qof$) to say that any automorphism is
the identity.) This theory is consistent as the reals can be expanded to
form a model of any finite subset of the theory. But a model of the
theory must have cardinality at least $\gk$. (Since we do not have the
downward L\"owenheim-Skolem theorem, we cannot assert that there is a
model of cardinality $\gk$.)

	In \cite{IV} and \cite{Sh84}, the compactness of 
two interesting
logics is established under certain set-theoretic hypotheses.  The
logics are those obtained from first order logic by adding quantifiers
which range over automorphisms of definable Boolean algebras or which
range over automorphisms of definable ordered fields. Instead of the
weaker version of dealing with automorphisms, it is also possible to
deal with a quantifier which says that two Boolean algebras are
isomorphic or that two ordered fields are isomorphic. The key step in
proving these results lies in establishing the following theorems. (By
definable we shall mean definable with parameters).

\begin{theorem}
\label{dmd-thm}
Suppose $\gl$ is a regular cardinal and both $\dmd(\gl)$ and
$\dmd(\set{\ga < \gl^+}{\cf \ga = \gl})$ hold. Then if $T$ is any
consistent theory and $|T| < \gl$, there is a model $M$ of $T$ of cardinality $\gl^+$ with
the following properties:
\begin{trivlist}
\item[(i)] If $B$ is a Boolean algebra definable 
 in $M$, then every automorphism of $B$ is definable.

\item[(ii)] $M$ is $\gl$-saturated.

\item[(iii)] Every non-algebraic type of cardinality $< \gl$ is
realized in $M$ by $\gl^+$ elements. 
\end{trivlist}
\end{theorem}

\begin{theorem}
\label{dmdfield}
Suppose $\gl$ is a regular cardinal and both $\dmd(\gl)$ and
$\dmd(\set{\ga < \gl^+}{\cf \ga = \gl})$ hold. Then if $T$ is any
consistent theory and $|T| < \gl$, there is a model $M$ of $T$ of cardinality $\gl^+$ with
the following properties:
\begin{trivlist}
\item[(i)] If $F$ is an ordered field definable
 in $M$ then every automorphism of $F$ is definable  and every
isomorphism between definable ordered fields is definable.

\item[(ii)] $M$ is $\gl$-saturated.

\item[(iii)] Every non-algebraic type of cardinality $< \gl$ is
realized in $M$ by $\gl^+$ elements. 

\item[(iv)] Every definable dense linear order is not the union of
$\gl$ nowhere dense sets. 
\end{trivlist}
\end{theorem}

	These theorems are 
proved in \cite{IV} (in
\cite{IV} section 9, Theorem~\ref{dmd-thm} is proved from GCH)  
although there is not an explicit statement of them there
. In order to show the
desired compactness result (from the assumption that there are
unboundedly many cardinals $\gl$ as in the theorem statements) it is
enough to use (i). However in our work on Boolean algebras we will
need the more exact information above. Let us notice how the
compactness of the various languages follow from these results. Since
the idea is the same in all cases just consider the case of Boolean
algebras. 

	First we will describe the logic ${\rm L}(\Qba)$. We  add
second order variables (to range over automorphisms of Boolean
algebras) and a quantifier $\Qba$ whose intended interpretation is
that there is an automorphism of the Boolean algebra.
More formally if $\gTh(f)$, $\gf(x), \psi(x,y), \gr(x, y)$ are
formulas (where $f$ is a second order variable, $x$ and $y$ are first
order variables and the formulas may have other variables) then
$$\Qba f\,(\gf(x), \psi(x,y))\,\,\gTh(f)$$ 
is a formula. In a model $M$
the tuple $(\gf(x), \psi(x,y))$ {\em defines} a Boolean algebra (where
parameters from $M$ replace the hidden free variables of $\gf(x),
\psi(x,y)$) if $\psi(x, y)$ defines a partial order $<$ on $B = \set{a
\in M}{M\models \gf[a]}$ so that $(B; <)$ is a Boolean algebra. A
model $M$ satisfies the formula $\Qba f\,(\gf(x),
\psi(x,y))\,\,\gTh(f)$ (where parameters from $M$ have been
substituted for the free variables), if whenever $(\gf(x), \psi(x,
y))$ define a Boolean algebra $(B; <)$ then there is an automorphism
$f$ of $(B; <)$ such that $M\models\gTh (f)$. (It is easy to extend
the treatment to look at Boolean algebras which are definable on
equivalence classes, but we will avoid the extra complication.) We can
give a more colloquial description of the quantifier $\Qba$ by saying
the interpretation of $\Qba$ is that ``$\Qba f (B) \ldots$'' holds if
there is an automorphism of the Boolean algebra $B$ so that \ldots.
We will describe some of the other logics we deal with in this looser
manner. For example, we will want to consider the quantifier $\Qof$
where $\Qof f (F_1, F_2) \ldots$ holds if there is is an isomorphism
$f$ from the ordered field $F_1$ to the ordered field $F_2$ such that
\ldots.

        The proof of compactness for ${\rm L}(\Qba)$ follows easily from
theorem 1.1. By expanding the language\footnote{More exactly,
 for every model $M$ define a model $M^*$ with universe
$M\cup\{f: f$ a
partial function from $\vert M\vert$ to $\vert M\vert\}$ with the relations
of $M$ and the unary predicate $P$,
$P^{M'}=\vert M\vert$,
and the ternary predicate $R$,
$R=\{(f, a, b): f\in {}^MM$,
$a\in M, b=f(a)\}$.
We shall similarly transform a theory $T$ to $T'$ and consider
automorphism only of structures $\subseteq P$.}
we can assume that there is a ternary relation
$R(\und,\und,\und)$ so that the theory says that any first order
definable function is definable by $R$ and one parameter.  By the
ordinary compactness theorem if we are given a consistent theory in this
logic then there is a model of the theory where all the sentences of the
theory hold if we replace automorphisms by definable automorphsisms in
the interpretation of $\Qba$, since quantification over definable
automorphisms can be replaced by first order quantification. Then we can
apply the theorem to get a new model elementarily equivalent to the one
given by the compactness theorem in which definable automorphisms and
automorphisms are the same.

In the following we will make two assumptions about all our theories.
First that all definable partial functions are in fact defined by a
fixed formula (by varying the parameters). Second we will always assume
that the language is countable except for the constant symbols.

	In this paper we will attempt to get compactness results
without recourse to $\dmd$, i.e., all our results will be in ZFC. We
will get the full result for the language where we quantify over
automorphisms (isomorphisms) of ordered fields in Theorem 6.4.
  Unfortunately we
are not able to show that the language with quantification over automorphisms
of Boolean algebras is compact, but will have to settle for a close
relative of that logic. This is theorem 5.1.
In section 4 we prove we can construct models in which all relevant
automorphism are somewhat definable:
4.1, 4.8 for BA, 4.13 for ordered fields.

        The reader may wonder why these results are being proved now,
about 10 years after the results that preceeded them. The key
technical innovation that has made these results possible is the
discovery of $\dmd$-like principles which are true in ZFC. These
principles, which go under the common name of the Black Box, allow one
to prove, with greater effort, many of the results which were
previously known to follow from $\dmd$ (see the discussion in
\cite{[Sh-e]} for more details). There have been previous applications of
the Black Box to abelian groups, modules and Boolean Algebras ---
often building objects with specified endomorphism rings. This
application goes deeper both in the sense that the proof is more
involved and in the sense that the result is more surprising. The
investigation is continued in \cite{[Sh 384]}, \cite{[Sh 482]}.

	In this paper we will also give a new proof of the compactness
of another logic --- the one which is obtained when a quantifier
$Q_{{\Brch}}$ is
added to first order logic which says that a level tree (definitions
will be given later) has an infinite branch. This logic was previously
shown to be compact --- in fact it was the first logic shown in ZFC to
be compact which is stronger than first order logic on countable
structures --- but our proof will yield a somewhat stronger result and
provide a nice illustration of one of our methods. (The first logic
stronger than first order logic which was shown to be compact was the
logic which expresses that a linear order has cofinality greater than
$\go_1$ \cite{cof}.) This logic, $L(Q_{{\Brch}})$,
 has been used by Fuchs-Shelah
\cite{316} to prove the existence of nonstandard uniserial modules
over (some) valuation domains. The proof uses the compactness of the tree
logic to transfer results proved using $\dmd$ to ZFC results.  Eklof
\cite{Ek} has given an explicit version of this transfer method and was
able to show that it settles other questions which had been raised.
(Osofsky \cite{Os1}, \cite{Os2} has found ZFC constructions which avoid
using the model theory.)

Theorem's 3.1 and 3.2 contain parallel results for Boolean algebra's and
fields.  They assert the existence of a theory (of sets) $T_1$ such that
in each model$M_1$ of $T_1$, $P(M_1)$ is a model $M$ of the first order
theory $T$ such that for every Boolean Algebra (respectively field)
defined in $M$, every automorphism of the Boolean algebra (repectively
field) that is definable in $M_1$ is definable in $M$.  Moreover, each
such $M_1$ has an elementary extension one of whose elements is a
pseudofinite set $a$ with the universe of $M_1$ contained in $a$ and with
$t(a/M_1)$ is definable over the empty set.  This result depends on the
earlier proof of our main result assuming $\diamondsuit$ and
absoluteness.  Theorem 4.1 uses the Black Box to construct a model $\GC$
of $T_1$ so that for any automorphism $f$ of a Boolean algebra $B =
P(\GC)$ there is a pseudofinite set $c$ such that for any atom $b \in B$,
$f(b)$ is definable from $b$ and $c$.  Theorem 4.13 is an analogous but
stronger result for fields showing that for any $b$, $f(b)$ is definable
from $b$ and $c$.  In Lemma 4.7, this pointwise definability is extended
by constructing a pseudo-finite partition of atoms of the Boolean algebra
(respectively the elements of the field) such that $f$ is definable on
each member of the partition.  In Theorem 5.1 for Boolean algebras and
6.4 for fields this local definability is extended to global
definability.

\subsection{Outline of Proof}

We want to build a model $A$ of a consistent $L(\Qof)$ theory $T$ which
has only definable isomorphisms between definable ordered fields. By the
ordinary compactness theorem, there is a non-standard model $\GC$ of an
expansion of a weak set theory (say ${\rm ZFC}^-$) which satisfies that
there is a model $A_1$ of $T$. So $A_1$ would be a model of $T$ if the
interpretation of the quantifier $\Qof$ were taken to range over
isomorphisms which are internal to $\GC$.  We can arrange that $A_1$
will be the domain of a unary predicate $P$. Then our goal is to build
our non-standard model $\GC$ of weak set theory in such a way that every
external isomorphism between definable ordered subfields of $P(\GC)$ is
internal, i.e., definable in $\GC$.

	The construction of $\GC$ is a typical construction with a
prediction principle, in this case the Black Box, where we kill
isomorphisms which are not pointwise definable over a set which is
internally finite (or synonymously, {\em pseudofinite}). A predicted
isomorphism is killed by adding an element which has no suitable
candidate for its image.  One common problem that is faced in such
constructions is the question ``how do we ensure no possible image of
such an element exists?''. To do this we need to omit some types. Much
is known about omitting a type of size $\gl$ in models of power $\gl$
and even $\gl^+$. But if say $2^{\gl} > \gl^{++}$, we cannot omit a
dense set of types of power $\gl$. So without instances of GCH we are
reduced to omitting small types, which is much harder. 
To omit the small types we will  use techniques which originated in
``Classification Theory''.  In the construction we will have for some
cardinal $\gth$ the type of any element does not split over a set of
cardinality less than $\gth$ (see precise definitions below). This is
analogous to saying the model is $\gth$-stable (of course we are
working in a very non-stable context).  The element we will add will
have the property that its image (if one existed) would split over
every set of cardinality $< \gth$.

	The final problem is to go from pointwise definablity to
definability. The first ingredient is a general fact about
$\ha_0$-saturated models of set theory. We will show for any
isomorphism $f$ that there is a large (internal) set $A$ and a
pseudofinite sequence of one-one functions $(f_i \colon i < k^*)$
which cover $f \rest A$ in the sense that for every $a \in A$ there is
$i$ so that $f_i(a) = f(a)$. Using this sequence of functions it is
then possible to define $f$ on a large subset of $A$. Finally, using
the algebraic structure, the
definition extends to the entire ordered field.

 	In this paper we will need to use the following principle. In
order to have the cleanest possible statement of our results (and to
conform to the notation in \cite{[Sh-e]}), we will state our results
using slightly non-standard notation. To obtain the structure ${\rm
H}_\chi(\gl)$, we first begin with a set of ordered urelements of
order type $\gl$ and then form the least set containing each urelement
and closed under formation of sets of size less than $\chi$. In a
context where we refer to ${\rm H}_\chi(\gl)$ by $\gl$ we will mean
the urelements and not the ordinals. In practice we believe that 
in a given context there will be no confusion.

\begin{theorem}
Suppose $\gl = \mu^+$, $\mu = \gk^\gth = 2^\gk$, $\chi$ is a regular
cardinal, $\gk$ is a strong limit cardinal, $\gth < \chi < \gk$, $\gk
> \cf \gk = \gth \geq \ha_0$ and $S \se \set{\gd < \gl}{\cf \gd =
\gth }$ is stationary. Let $\rho$ be some cardinal greater than $\gl$.
Then we can find $W = \set{(\bar{M}^\ga, \eta^\ga)}{\ga < \ga(\ast)}$
(actually a sequence)\footnote{ in our
case $\ga(*)=\gl$ is fine}, a function $\gz:\ga(\ast) \to S$ and $(C_\gd
\colon \gd \in S)$ such that:

\begin{trivlist}
\item[(a1)] $\bar{M}^\ga = (M^\ga_i \colon i \leq \gth)$ is an
increasing continuous elementary chain, each $M_i^\ga$ is a model
belonging to ${\rm H}_\chi(\gl)$ (and so necessarily has cardinality
less than $
\chi$), $M^\ga_i \cap \chi$ is an ordinal, $\eta^\ga \in {}^\gth\!\gl$
is increasing with limit $\gz(\ga) \in S$, for $i<\theta$,
$\eta^\ga\rest i \in
M^\ga_{i+1}$, $M^\ga_i \in {\rm H}_\chi(\eta^\ga(i))$ and $(M^\ga_j
\colon j \leq i) \in M_{i+1}^\ga$.

\item[(a2)] For any set $X \se \gl$ there is $\ga$ so that $M^\ga_\gth
\equiv_{\gl \cap M^\ga_\gth} ({\rm H}(\rho), \in, <, X)$, where $<$ is a
well ordering of ${\rm H}(\rho)$ and $M\equiv_A N$
means $(M, a)_{a\in A},
(N, a)_{a\in A}$ are elementarily equivalent.

\item[(b0)] If $\ga \neq \gb$ then $\eta^\ga \neq \eta^\gb$.

\item[(b1)] If $\set{\eta^\ga \rest i}{i < \gth} \se M^\gb_\gth$ and
$\ga \neq \gb$ then $\gz(\ga) < \gz(\gb)$.

\item[(b2)] If $\eta^\ga \rest (j + 1) \in M^\gb_\gth$ then $M^\ga_j \in
M^\gb_\gth$. 

\item[(c2)] $\bar{C} = (C_\gd \colon \gd \in S)$ is such that each
$C_\gd$ is a club subset of $\gd$ of the order type $\gth$.

\item[(c3)] Let $C_\gd = \set{\grg_{\gd,i}}{i < \gth}$ be an increasing
enumeration. For each  $\ga < \ga(\ast)$ there is $( (\grg^-_{\ga, i},
\grg^+_{\ga,i}) \colon i < \gth)$ such that:
$\grg^-_{\ga,i} \in M^\ga_{i+1}$, 
$M^\ga_{i+1} \cap \gl \se
\grg^+_{\ga,i}$, $\grg_{\gz(\ga),i} < \grg^-_{\ga,i} < \grg^+_{\ga,i} <
\grg_{\gz(\ga),i+1}$ and if $\gz(\ga) = \gz(\gb)$ and $\ga \ne \gb$ then for
every large enough $i < \gth$, 
 $[\grg^-_{\ga,i}, \grg^+_{\ga,i}) \cap [\grg^-_{\gb,i},
\grg^+_{\gb,i}) = \emptyset$. Furthermore for all  $i$, the
sequence $(\grg^-_{\ga, j} \colon j  < i)$ is in $M^\ga_{i+1}$.

\end{trivlist}
\end{theorem}

	This principle, which is one of the Black Box principles is a
form of $\dmd$ which is a theorem of ZFC. This particular principle is
proved in \cite{[Sh-e]}III 6.13(2). The numbering here is chosen to
correspond with the numbering there. Roughly speaking clauses (a1) and
(a2) say that there is a family of elementary substructures which
predict every subset of $\gl$ as it sits in $\mbox{\rm H}(\gl)$. (We
will freely talk about a countable elementary substructure {\em
predicting} isomorphisms and the like.) The existence of such a family
would be trivial if we allowed all elementary substructures of
cardinality less than $\chi$. The rest of the clauses say that the
structures are sufficiently disjoint that we can use the information
that they provide without (too much) conflict.

	The reader who wants to follow the main line of the arguments
without getting involved (initially) in the complexities of the Black
Box can substitute $\dmd(\gl)$ for the Black Box. Our proof of the
compactness of ${\rm L}(\Qof)$ does not depend on
Theorem~\ref{dmdfield}, so even this simplification gives a new proof
of the consistency of the compactness of  ${\rm L}(\Qof)$. Our work on
Boolean algebras does require Theorem~\ref{dmd-thm}.

	The results in this paper were obtained while the first author
was visiting the Hebrew University in Jerusalem. He wishes to thank
the Institute of Mathematics for its hospitality.

\section{Non-splitting extensions}
\setcounter{theorem}{0}

        In this section $\gth$ will be a fixed regular cardinal.
Our treatment is self contained but the reader can look at
\cite{[Sh-a]}.

 \medskip

\noindent {\sc Definition}. If $M$ is a model and $X, Y, Z \se M$, then $X/Y$
{\em does not split over} $Z$ if and only if for every finite $d \se Y$
the type of $X$ over $d$ (denoted either tp($X/d$) or $X/d$) depends
only on the type of $d$ over $Z$. \medskip

	We will use two constructions to guarantee that types will not
split over small sets. The first is obvious by definition. (The type
of $A/B$ is {\em definable over} $C$ if for any tuple $\bar{a} \in A$
and formula $\gf(\bar{x}, \bar{y})$ there is a formula
$\psi(\bar{y})$ with parameters from $C$ so that for any $\bar{b}$,
$\gf(\bar{a}, \bar{b})$ if and only if $\psi(\bar{b})$.)

\begin{proposition} If $X/Y$ is definable over $Z$ then $X/Y$ does not split
over $Z$. 
\end{proposition}

\noindent
{\sc Definition}. Suppose $M$ is a model and $X, Y \se M$. Let $D$ be an
ultrafilter on $X^\ga$. Then the $\mbox{Av}(X,D,Y)$ (read the
($D$-)average type that $X^\ga$ realizes over $Y$) is the type $p$ over
$Y$ defined by: for $\bar{y}\subseteq Y,  \gf(\bar{z}, \bar{y}) \in p$
if and only if $\set{\bar{x} \in X^\ga}{\gf(\bar{x}, \bar{y}) \mbox{
holds}} \in D$. We will omit $Y$ if it is clear from context. Similarly
we will omit $\ga$ and the ``bar'' for singletons, i.e., the case $\ga
=1$.\medskip

The following two propositions are clear from the definitions. 

\begin{proposition}
\label{ult}
 If $\bar{a}$ realizes $\mbox{Av}(X,D,Y)$ then
$\bar{a}/Y$ does not split over $X$. Also if there is $Z$ such that for
$\bar{b} \in X$ the type of $\bar{b}/Y$ does not split over $Z$, then
$\mbox{Av}(X,D,Y)$ \dns\ $Z$.
\end{proposition}

\begin{proposition}
\label{trans}
\begin{trivlist}
\item[(i)] Suppose $A/B$ \dns\ $D$, $B \se C$ and  $C/B\cup A$ \dns\
$D\cup A$ then $A\cup C/B$
\dns\ $D$.

\item[(ii)] Suppose that $(A_i \colon i < \gd)$ is an increasing chain
and for all $i$, $A_i/B$ \dns\ $C$ then $\bigcup_{i < \gd} A_i/B$ \dns\
C.
\item[(iii)] $X/Y$ \dns\ $Z$ \iff\ $X/{\rm dcl}(Y \cup Z)$ \dns\ $Z$.
Here ${\rm dcl} (Y \cup Z)$ denotes the definable closure of $Y \cup Z$.

\item[(iv)] If $X/Y$ \dns\ $Z$ and $Z \se W$, then $X/Y$ \dns\ $W$
\end{trivlist}
\end{proposition}

\noindent {\sc Definition}. Suppose $M_1 \prec M_2$ are models. Define $M_1
\slt M_2$, if for every $X \se M_2$ of cardinality less than $\gth$
there is $Y \se M_1$ of cardinality less than $\gth$ so that $X/M_1$
\dns\ $Y$. (If $\gth$ is regular, then we only need to consider the
case where $X$ is finite.)

\begin{proposition} 
Assume that $\theta$ is a regular cardinal (needed for (2) only) and
that all models are models of some fixed theory with Skolem functions
(although this is needed for (3) only).
\begin{trivlist}
\item[1.] $\slt$ is transitive and for all $M$, $M \slt M$.
\item[2.] If $(M_i \colon i < \gd)$ is a $\slt$-increasing chain, then
for all $i$ $M_i \slt \bigcup_{j < \gd} M_j$.
\item[3.] Suppose $M_2$ is generated by $M_1 \cup N_2$ and $N_1 = M_1
\cap N_2$. (Recall that we have Skolem functions.) If $|N_1| < \gth$
and $N_2/M_1$ \dns\ $N_1$, then $M_1 \slt M_2$.
\end{trivlist}
\end{proposition}

An immediate consequence of these propositions is the following
proposition.

\begin{proposition}
\label{saturated}
Suppose $M \slt N$, then there is a $\gth$-saturated model $M_1$ such
that $N \prec M_1$ and $M \slt M_1$.
\end{proposition}

\proof By the lemmas it is enough to show that given a set $X$ of
cardinality $< \gth$ and a type $p$ over $X$ we can find a realization $a$ of
that type so that $a/N$ \dns\ $X$. Since we have Skolem functions
every finite subset of $p$ is realized by an element whose type over
$N$ does not split over $X$, namely an element  of the Skolem hull of
$X$. So we can take $a$ to realize an average type of these elements. \fin

\section{Building New Theories}
\setcounter{theorem}{0}

	The models we will eventually build will be particular
non-standard models of an enriched version of ${\rm ZFC}^-$.  (Recall
${\rm ZFC}^-$ is ZFC without the power set axiom and is true in the sets
of hereditary cardinality $< \gk$ for any regular uncountable cardinal
$\gk$.) The following two theorems state that appropriate theories
exist. 

\begin{theorem}
\label{bathy}
Suppose $T$ is a theory in a language which is countable except for
constant symbols, $P_0$ a unary predicate
 so that in every model $M$ of $T$ every definable
automorphism of a definable atomic Boolean algebra$\,\se P^M_0$
 is definable by a fixed
formula (together with some parameters). Then there is $T_1$ an
expansion of ${\rm ZFC}^-$ in a language which is countable except for
constant symbols with a unary predicate $P_0$ so that if $M_1$
is a model of $T_1$ then  $P_0(M_1)$ is  a model $M$ of $T$
(when restricted to the right vocabulary) and the
following are satisfied to be true in $M_1$.
\begin{trivlist}
\item[(i)] Any automorphism of a  definable (in $M$) atomic Boolean
algebra contained in  $P_0(M)$
 which is definable in $M_1$ is definable in $M$.
\item[(ii)] $M_1$ (which is a model of ${\rm ZFC}^-$) satisfies for
some regular cardinal $\mu$ (of $M_1$), $|M| = \mu^+$, $M$ is 
$\mu$-saturated and every non-algebraic type (in the language of $M$)
of cardinality $< \mu$ is realized in $M$ by $\mu^+$ elements of $M$.
\item[(iii)] $M \in M_1$.
\item[(iv)] $M_1$ satisfies the separation scheme for all formulas
(not just those of the language of set theory).
\item[(v)] $M_1$ has Skolem functions.
\item[(vi)] For any $M_1$ there is an elementary extension $N_1$ so that
the universe of $M_1$ is
 contained in $N_1$ in a pseudofinite set (i.e., one which is
finite in $N_1$) whose  type over the universe of $M_1$ is definable
over the empty set.
\end{trivlist}

\end{theorem}

\proof We  first consider a special case. Suppose that there is a
cardinal $\gl$ greater than the cardinality of $T$ satisfying the
hypothess of Theorem~\ref{dmd-thm} (i.e., both $\dmd(\gl)$ and
$\dmd(\set{\ga < \gl^+}{\cf \ga = \gl})$ hold). Then we could choose
$\gk$ a regular cardinal greater than $\gl^+$. Our model $M_1$ will be
taken to be a suitable expansion of $\her{\gk}$ where the interpretation
of the unary predicate $P_0$ is the model $M$ guaranteed by
Theorem~\ref{dmd-thm} and $\mu = \gl$, $\mu^+ = \gl^+$. Since any
formula in the enriched language is equivalent in $M_1$ to a formula of
set theory together with parameters from $M_1$, $M_1$ will also
satisfy (iv).

What remains is to ensure (v) and that appropriate elementary extensions
always exist i.e. clause (vi).  To achieve this we will expand the
language by induction on $n$. Let $L_0$ be the language consisting of
the language of $T$, $\{P\}$, and Skolem functions and $M_1^0$ be any
expansion by Skolem functions of the structure on $\her{\gk}$ described
above. Fix an index set $I$ and an ultrafilter $D$ on $I$ so that there
is $a \in \her{\gk}^I/D$ such that $\her{\gk}^I/D
\models
\mbox{`` $a$ is finite''}$ and for all $b \in \her{\gk}$,
$\her{\gk}^I/D \models b
\in a$. Let $N_1^0 = {M_1^0}^I/D$.

	  Then for every formula $\gf(y, x_1, \ldots, x_n)$ of $L_0$
which does not involve constants, add a new $n$-ary relation $R_\gf$.
Let $L_{0.5}$ by the language containing all the $R_\gf$.  Let
$M_1^{0.5}$ be the $L_{0.5}$ structure with universe $H(\kappa)$
obtained by letting for all $b_1, \ldots, b_n$, $M_1^{0.5}\models
R_\gf[b_1, \ldots, b_n]$ if and only if $N_1^0
\models \gf[a, b_1, \ldots, b_n]$. Let $L_1$ be an extension of
$L_{0.5}$ by Skolem functions and let $M_1^1$ be an expansion of
$M_1^{0.5}$ by Skolem functions. Condition (iv)
still holds as it holds for any expansion of
$(H(\kappa), \in)$. We now let $N^1_1 = {M^1_1}^I/D$ and
continue as before. 

Let $L = \bigcup_{n <\go} L_n$ and $M_1 = \bigcup_{n<\go} M_1^n$ (i.e.
the least common expansion; the universe stays the same).  Let $T_1$ be
the theory of $M_1$. As we have already argued $T_1$ has properties
(i)--(iv). It remains to see that any model of $T_1$ has the desired
extension property. First we consider $M_1$ and let $N_1 = {M_1}^I/D$.
Then the type of $a$ over $M_1$ is definable over the empty set using
the relations $R_\gf$ which we have added.  Since $T_1$ has Skolem
functions, for any model $A_1$ of $T_1$ there will be an extension $B_1$
of $A_1$ generated by $A_1$ and an element realizing the definable type
over $A_1$.

In the general case, where we may not have the necessary hypotheses of
Theorem~\ref{dmd-thm}, we can force with a notion of forcing which adds
no new subsets of $|T|$ to get some $\gl$ satisfying hypothess of
Theorem~\ref{dmd-thm} (alternately, we can use L[$A$] where $A$ is a
large enough set of ordinals). Since the desired theory will exist in an
extension it already must (as it can be coded by a subset of $|T|$)
exist in the ground model. \fin

Later on we will be juggling many different models of set theory. The
ones which are given by the Black Box, and the non-standard ones which
are models of $T_1$. When we want to refer to notions in models of
$T_1$, we will use words like ``pseudofinite'' to refer to sets which
are satisfied to be finite in the model of $T_1$.

In the same way as we proved the last theorem we can show the following
theorem.

\begin{theorem}
\label{ofthy}
Suppose $T$ is a theory in a language with a unary predicate $P_0$ and
which is countable except for constant symbols so that in every model
$M$ of $T$ every definable isomomrphism between definable ordered
fields$\subseteq P_0^M$ is definable by a fixed formula (together with
some parameters). Then there is $T_1$ an expansion of ${\rm ZFC}^-$ in a
language which is countable except for constant symbols with a unary
predicate $P_0$ so that if $M_1$ is a model of $T_1$ then $P_0(M_1)$ is a
model $M$ of $T$ and the following are satisfied to be true in $M_1$.
\begin{trivlist}
\item[(i)] Any isomorphism of definable (in $M$) ordered fields
contained in $P_0$ which is definable in $M_1$ is definable in $M$.
\item[(ii)] $M \in M_1$.
\item[(iii)] $M_1$ satisfies the separation scheme for all formulas
(not just those of the language of set theory).
\item[(iv)] $M_1$ has Skolem functions.
\item[(v)] For any $M_1$ there is an elementary extension $N_1$ so that
the universe of $M_1$ is contained in $N_1$ in a pseudofinite set (i.e.,
one which is finite in $N_1$) whose type over the universe of $M_1$ is
definable over the empty set.
\end{trivlist}

\end{theorem}

Since we have no internal saturation conditions this theorem can be
proved without recourse to Theorem~\ref{dmdfield} (see the next section
for an example of a similar construction).

\subsection{A Digression}
The method of expanding the language to get extensions which realize a
definable type is quite powerful in itself. We can use the method to
give a new proof of the compactness of a logic which extends first order
logic and is stronger even for countable structures. This subsection is
not needed in the rest of the paper.

\begin{lemma}
\label{tree-lemma}
Suppose that $N$ is a model.  Then there is a consistent expansion of
$N$ to a model of a theory $T_1$ with Skolem functions so that for every
model $M$ of $T_1$ there is $a$ so that $a/M$ is definable over the
empty set and in $M(a)$ (the model generated by $M$ and $\{a\}$) for
every definable directed partial ordering $<$ of $M$ without the last
element there is an element greater than any element of $M$ in the
domain of $<$.  Furthermore the cardinality of the language of $T_1$ is
no greater than that of $N$ plus $\ha_0$.
\end{lemma}

\proof Fix a model $N$. Choose $\gk$ and an ultrafilter $D$ so
that for every directed partial ordering $<$ of $N$ without the last
element there is an element of $N^\gk/D$ which is greater than every
element of $N$ (e.g., let $\gk = |N|$ and $D$ be any regular ultrafilter
on $\gk$). Fix an element $a \in N^\gk/D \setminus N$. Abusing notation
we will let $a :\gk\to N$ be a function representing the element $a$.
The new language is defined by induction on $\go$. Let $N = N_0$. There
are three tasks so we divide the construction of $N_{n+1}$ into three
cases. If $n \equiv 0 \bmod 3$, expand $N_n$ to $N_{n+1}$ by adding
Skolem functions. If $n \equiv 1 \bmod 3$, add a $k$-ary relation
$R_\gf$ for every formula of arity $k+1$ and let $R_\gf(b_0, \ldots,
b_{k-1})$ hold if and only if $\gf(b_0, \ldots, b_{k-1}, a)$ holds in
${N_n}^\gk/D$.

If $n \equiv 2 \bmod 3$, we ensure that there is an upper bound to every
definable directed partial orders without last element in $N_n$. For
each $k+2$-ary formula $\gf(x_0, \ldots, x_{k-1}, y, z)$ we will add a
$k+1$-ary function $f_\gf$ so that for all $\bar{b}$ if $\gf(\bar{b},
y,z)$ defines a directed partial order without last element then
$f_\gf(\bar{b},a)$ is greater in that partial order than any element of
$N$. Notice that there is something to do here since we must define
$f_\gf$ on $N$ and then extend to $N^\gk/D$ using the ultraproduct. For
each such $\bar{b}$ choose a function $c:\gk\to N$ so that $c/D$ is an
upper bound (in the partial order) to all the elements of $N$. Now
choose $f_\gf$ so that $f_\gf(\bar{b}, a(i)) = c(i)$. Let $T_1$ be the
theory of the expanded model.

Suppose now that $M$ is a model of $T_1$. The type we want is the type
$p$ defined by $\gf(c_1, \ldots, c_{k-1},x) \in p$ if and only if $M
\models R_\gf(c_1, \ldots, c_n)$. \fin

In the process of building the new theory there are some choices made of
the language. But these choices can be made uniformly for all models. We
will in the sequel assume that such a uniform choice has been made.

>From this lemma we can prove a stronger version of a theorem from
\cite{II}, which says that the language which allows quantification over
branches of level-trees is compact. A {\em tree} is a partial order in
which the predecessors of any element are totally ordered. A {\em
level-tree} is a tree together with a ranking function to a directed
set. More exactly a {\em level-tree} is a model $(A: U, V, <_1, <_2, R)$
where:
\begin{enumerate}
\item $A$ is the union of $U$ and $V$;
\item $<_1$ is a partial order of $U$ such that for every $u \in U$
$\set{y \in U}{y <_1 u}$ is totally ordered by $<_1$;
\item $<_2$ is a directed partial order on $V$ with no last element; 
\item $R$ is a function from $U$ to $V$ which is strictly order preserving.
\end{enumerate}
The definiton here is slightly more general than in \cite{II}. In
\cite{II} the levels were required to be linearly ordered. Also what we
have called a ``level-tree'' is called a ``tree'' in \cite{II}. A {\em
branch} $b$ of a level-tree is a maximal linearly-ordered subset of $U$
such that $\set{R(t)}{t \in b}$ is unbounded in $V$. We will refer to
$U$ as the tree and $V$ as the levels. For $t \in U$ the {\em level} of
$t$ is $R(t)$.

A tuple of formulas (which may use parameters from $M$),
$$(\gf_1(x),\gf_2(x), \psi_1(x, y), \psi_2(x, y), \rho(x,y, z)),$$
defines a level-tree in a model $M$ if $$(\gf_1(M) \cup
\gf_2(M);\gf_1(M),\gf_2(M), \psi_1(x, y)^M, \psi_2(x, y)^M,
\rho(x,y,z)^M)$$
is a level tree. (There is no difficulty in extending the treatment to
all level-trees which are definable using equivalence relations.)

Given the definition of a level-tree, we now define an extension of
first order logic by adding second order variables (to range over
branches of level-trees) and a quantifier $Q_{\rm Brch}$, such that
$$\Qb \,b (\gf_1(x),\gf_2(x), \psi_1(x, y), \psi_2(x, y), \rho(x,y, z))
\gTh(b)$$
says that if $(\gf_1(x),\gf_2(x), \psi_1(x, y), \psi_2(x, y), \rho(x,y,
z))$ defines a level-tree then there is a branch $b$ of the level-tree
such that $\gTh(b)$ holds. 

In \cite{II}, it is shown that (a first-order version of) this logic was
compact.  This is the first language to be shown (in ZFC) to be fully
compact and stronger than first order logic for countable structures.
In \cite{II}, the models are obtained at successors of regular
cardinals.

\begin{theorem}
The logic $L(Q_{\rm Brch})$ is compact. Furthermore every consistent
theory $T$ has a model in all uncountable cardinals $\gk > |T|$.
\end{theorem}

\proof By expanding the language we can assume that any model of $T$
admits elimination of quantifiers. (I.e., add new relations for each
formula and the appropriate defining axioms.)

For each finite $S \se T$, choose a model $M_S$ of $S$. For each $S$
choose a cardinal $\mu$ so that $M_S \in {\rm H}(\mu)$ and let $N_S$ be
the model $({\rm H}(\mu^+), M_S, \in)$, where $M_S$ is the
interpretation of a new unary predicate $P$ and the language of $N_S$
includes the language of $T$ (with the correct restriction to $M_S$). By
expanding the structure $N_S$ we can assume that the theory $T_S$ of
$N_S$ satisfies the conclusion of Lemma~\ref{tree-lemma}.  Furthermore
we note that $T_S$ satisfies two additional properties. If a formula
defines a branch in a definable level-tree contained in the domain of
$P$ then this branch is an element of the model $N_S$. As well $N_S$
satisfies that $P(N_S)$ is a model of $S$.

Now let $D$ be an ultrafilter on the finite subsets of $T$ such that for
all finite $S$, $\set{S_1}{S \se S_1} \in D$. Finally let $T_1$ be the
(first-order) theory of $\prod N_S/D$. If $N$ is any model of $T_1$ then
for any sentence $\gf \in T$, $N$ satisfies ``$P(N)$ satisfies $\gf$''.
If we can arrange that the only branches of an $L(Q_{\rm
Brch})$-definable (in the language of $T$) level-tree of $P(N)$ are
first order definable in $N$ then $P(N)$ satisfaction of an $L(Q_{\rm
Brch})$-formula will be the same in $N$ and the real world. Before
constructing this model let us note that our task is a bit easier than
it might seem.

\begin{claim}
Suppose that $N$ is a model of $T_1$ and every branch of an first-order
definable (in the language of $T$) level-tree of $P(N)$ is first order
definable in $N$, then every branch of an $L(Q_{\rm Brch})$-definable
(in the language of $T$) level-tree of $P(N)$ is first order definable
in $N$.
\end{claim}

\proof (of the claim) Since we have quantifier elimination we can prove by
induction on construction of formulas that satisfaction is the same in
$N$ and the real world and so the quantifier-elimination holds for
$P(N)$. In other words any $L(Q_{\rm Brch})$-definable level-tree is
first order definable.

It remains to do the construction and prove that it works. To begin let
$N_0$ be any model of $T_1$ of cardinality $\gk$. Let $\mu
\leq \gk$ be any regular cardinal. We will construct an increasing
elementary chain of models $N_\ga$ for $\ga <
\mu$ by induction. At limit ordinals we will take unions. If $N_\ga$
has been defined, let $N_{\ga +1} = N_\ga(a_\ga)$, where $a_\ga$ is as
guaranteed by Lemma~\ref{tree-lemma}. Now let $N = \bigcup_{\ga<\gm}
N_\ga$.

\begin{subclaim}
Suppose $X$ is any subset of $N$ which is definable by parameters. Then
for all $\ga$, $X \cap N_\ga$ is definable in $N_\ga$.
\end{subclaim}	

\proof  (of the subclaim) Suppose not and let $\gb$ be the least ordinal
greater than $\ga$ so that $X \cap N_\gb$ is definable in $N_\gb$. Such
an ordinal must exist since for sufficiently large $\gb$ the parameters
necessary to define $X$ are in $N_\gb$. Similarly there is $\grg$ such
that $\gb = \grg + 1$. Since $N_\gb$ is the Skolem hull of $N_\grg \cup
\{a_\grg\}$, there is $\bar{b} \in N_\grg$ and a formula $\gf(x,
\bar{y}, z)$ so that $X \cap N_\gb$ is defined by $\gf(x, \bar{b},
a_\grg)$. But by the definability of the type of $a_\grg$ over $N_\grg$
there is a formula $\psi(x, \bar{y})$ so that for all $a, \bar{c} \in
N_\grg$, $N_\gb \models \psi(a, \bar{c}) \mbox{ if and only if } \gf(a,
\bar{c}, a_\grg)$. Hence $\psi(x, \bar{b})$ defines $X \cap N_\grg$ in
$N_\grg$. 

It remains to see that every branch of a definable level tree is
definable. Suppose $(A; U, V,<_1, <_2, R)$ is a definable level-tree.
Without loss of generality we can assume it is definable over the empty
set. Let $B$ be a branch.  For any $\ga < \gm$ there is $c \in V$ so
that for all $d \in V \cap N_\ga$, $d <_2 c$. Since the levels of $B$
are unbounded in $V$, $B \cap N_\ga$ is not cofinal in $B$. Hence there
is $b \in B$ so that $B \cap N_\ga \se \set{a \in U}{a <_1 b}$. By the
subclaim $B \cap N_\ga$ is definable in $N_\ga$.

Since $\gm$ has the uncountable cofinality, by Fodor's lemma, there is
$\ga< \gm$ so that for unboundedly many (and hence all) $\grg < \gm$, $B
\cap N_\grg$ is definable by a formula with parameters from $N_\ga$. 
Fix a formula $\gf(x)$ with parameters from $N_\ga$ which defines $B
\cap N_\ga$.  Then $\gf(x)$ defines $B$.  To see this consider any $\grg
< \ga$ and a formula $\psi(x)$ with parameters from $N_\ga$ which
defines $B \cap N_\grg$.  Since $N_\ga$ satisfies ``for all $x$,
$\gf(x)$ if and only if $\psi(x)$'' and $N_\ga \prec N_\grg$, $\gf(x)$
also defines $B \cap N_\grg$.  \fin

\noindent {\sc Remark}. The compactness result  above is optimal as far as
the cardinality of the model is concerned. Any countable level-tree has
a branch and so there is no countable model which is $L(Q_{\rm
Brch})$-equivalent to an Aronszajn tree. By the famous theorem of
Lindstrom \cite{Li}, this result is the best that can be obtained for
any logic, since any compact logic which is at least as powerful as the
first order logic and has countable models for all sentences is in fact
the first order logic. The existence of a compact logic such that every
consistent countable theory has a model in all uncountable cardinals was
first proved by Shelah \cite{cof}, who showed that the first order logic
is compact if we add a quantifier $Q^{\rm cf}$ which says of a linear
order that its cofinality is $\go$. (Lindstrom's theorem and the logic
$L(Q^{cf})$ are also discussed in \cite{Ba}). The logic $L(\Qb)$ has the
advantage that it is stronger than the first order logic even for
countable models \cite{II}.

Notice in the proof above in any definable level-tree the directed set
of the levels has cofinality $\gm$. Since we can obtain any uncountable
cofinality this is also the best possible result.  
Also in the theorem above  we can demand just $\kappa \ge |T| + \aleph_1$


\section{The Models}
\setcounter{theorem}{0}

For the purposes of this section let $T$ and $T_1$ be theories as
defined above in 3.1 or 3.2 (for Boolean algebras or ordered fields). In
this section we will build models of our theory $T_1$. The case of
Boolean algebras and the case of ordered fields are similar but there
are enough differences that they have to be treated separately. We shall
deal with Boolean algebras first.

We want to approximate the goal of having every automorphism of every
definable atomic Boolean algebra in the domain of $P$ be definable. In
this section, we will get that they are definable in a weak sense. In
order to spare ourselves some notational complications we will make a
simplifying assumption and prove a weaker result. It should be apparent
at the end how to prove the same result for every definable atomic
Boolean algebra in the domain of $P$.

\medskip
\noindent {\sc Assumption.} Assume $T$ is the theory of an atomic
Boolean algebra on $P$ with some additional structure.

\begin{theorem} 
\label{Ba-thm}
There is a model $\GC$ of $T_1$ so that if $B = P(\GC)$ and $f$ is any
automorphism of $B$ as a Boolean algebra then there is a pseudofinite
set $c$ so that for any atom $b \in B$, $f(b)$ is definable from $b$ and
elements of $c$.
\end{theorem}

\proof  We will use the notation from the Black Box. In particular we
will use an ordered set of urelements of order type $\gl$. We can assume
that $\mu$ is larger than the cardinality of the language (including the
constants). We shall build a chain of structures $(\GC_\vare
\colon \vare < \gl)$ such that the universe of $\GC_\vare$  will be
 an ordinal $<\gl$ and the universe of $\GC=\bigcup_{\vare<\gl}
\GC_\vare$ will be $\gl$ (we can specify in a definable way what the
universe of $\GC_\vare$ is e.g. $\mu(1+\ga)$ is o.k.). We choose
$\GC_\vare $ by induction on $\vare$. Let $B_\vare = P(\GC_\vare)$.  We
will view the $(\bar M^\ga, \eta^\alpha)\in W$ in the Black Box as
predicting a sequence of models of $T_1$ and an automorphism of the
Boolean algebra.  (See the following paragraphs for more details on what
we mean by predicting.)

The construction will be done so that if $\vare \not\in S$ and $\vare <
\gz$ then $\GC_\vare \slt \GC_\gz$. (We will make further demands
later.) The model $\GC$ will be $\bigcup_{\vare<\gl}
\GC_\vare$.  When we are done, if $f$ is an automorphism of the Boolean
algebra then we can choose $(\bar M^\ga, \eta^\ga)\in W$ to code, in a
definable way, $f$ and the sequence $(\GC_\vare \colon \vare < \gl)$.

The limit stages of the construction are determined. The successor stage
when $\vare \notin S$ is simple. We construct $\GC_{\vare+1}$ so that
$\GC_\vare\slt \GC_{\vare+1}$,  $\GC_{\vare+1}$ is $\gth$-saturated
and there is a pseudofinite set in $\GC_{\vare+1}$ which contains
$\GC_\vare$. By the construction of the theory $T_1$ there is $c$, a
pseudofinite set which contains $\GC_\vare$ such that the type of $c$
over $\GC_\vare$ is definable over the empty set. Hence $\GC_\vare \slt
\GC_\vare(c)$.  By Proposition~\ref{saturated} there is $\GC_{\vare+1}$
which is $\gth$-saturated so that $\GC_\vare(c) \slt \GC_{\vare+1}$.
Finally by transitivity (Proposition~\ref{trans}), $\GC_\vare\slt
\GC_{\vare+1}$.

The difficult case occurs when $\vare \in S$ rename $\vare$ by $\gd$.
Consider $\ga$ so that $\gz(\ga) = \gd$. We are interested mainly in
$\alpha$'s which satisfy:\\ 
$(*)$ \ \ \ The model $M_\gth^\ga$ ``thinks'' it is of the form
$(\her\rho,\in, <, X)$ and by our coding yields (or predicts) a sequence
of structures $(\mon{D}_\gn\colon \gn < \gl)$ and a function $f_\ga$
from $\mon{D} = \bigcup_{\gn < \gl} \mon{D}_\gn$ to itself. (Of course
all the urelements in $M_\gth^\ga \cap \gl$ will all have order type
less than $\gd$.)\\ 
At the moment we will only need to use the function
predicted by $M_\gth^\ga$.

We will say an {\em obstruction occurs at} $\ga$ if (($*$) holds and) we
can make the following choices. If possible choose $N_\ga \se \GC_\gd$
so that $N_\ga \in M_0^\ga$ of cardinality less than $\gth$ and a
sequence of atoms $(a^\ga_i \colon i \in C_\gd $) so that (naturally
$a^\ga_i\in M^\ga_{i+1})\, a^\ga_i/\GC_{\gamma^-_{\ga,i}}$ \dns\ $N_\ga$
and $a_i^\ga \in [{\gamma^-_{\ga,i}}, {\gamma^+_{\ga,i}})$ and
$f_\ga(a^\ga_i)$ is not definable over $a$ and parameters from
$\GC_{\gamma^-_{\ga, i}}$. At ordinals where an obstruction occurs we
will take action to stop $f_\ga$ from extending to an automorphism of
$B=B^{\GC}$.

Notice that $N_\ga$ is contained in $M^\ga_0\se M_\gth^\ga$.

Suppose an obstruction occurs at $\ga$. Let $X_\ga$ be the set of finite
joins of the $\set{a^\ga_i}{i \in C_\gd }$.  In the obvious way, $X_\ga$
can be identified with the set of finite subsets of
$Y_\ga=\set{a^\ga_i}{i \in C_\gd} $.  Fix $U_\ga$ an ultrafilter on
$X_\ga$ so that for all $x\in X_\ga$, $\set{y \in X_\ga}{x \se y} \in
U_\ga$. Now define by induction on $\mbox{Ob}(\gd) = \set{\ga}{\gz(\ga)
= \gd \mbox{ and an obstruction occurs at }\ga}$ an element $x_\ga$ so
that $x_\ga$ realizes the $U_\ga$ average type of $X_\ga$ over $\GC_\gd
\cup\set{x_\gb}{\gb \in \mbox{Ob}(\gd), \gb < \ga}$. Then $\GC_{\gd + 1}$ is
the Skolem hull of $\GC_\gd \cup \set{x_\ga}{\ga \in \mbox{Ob}(\gd)}$.

We now want to verify the inductive hypothesis and give a stronger
property which we will use later in the proof. The key is the following
claim.

\begin{claim}
Suppose $\ga_0, \ldots, \ga_{n-1} \in \Ob$, then for all but a bounded set
of $\grg < \gd$, $(\bigcup_{k < n} Y_{\ga_k})/\GC_\grg$ \dns\ $\bigcup_{k<n}
N_{\ga_k} \cup \bigcup_{k<n} (\GC_\grg \cap Y_{\ga_k})$.
\end{claim}

\proof (of the claim) Suppose $\grg$ is large enough so that for all $ m
\not=k < n$, $[{\gamma^-_{\ga_m,i}}, {\gamma^+_{\ga_m,i}}) \cap
[{\gamma^-_{\ga_k,i}}, {\gamma^+_{\ga_k,i}}) = \emptyset$, whenever
$\grg^-_{\ga_m, i}\geq \grg$ (recall clause (c3) of the Black Box). It
is enough to show by induction on $\grg \leq \gs < \gd$ that $(\bigcup_{k <
n} Y_{\ga_k}) \cap \GC_\gs$ has the desired property.  For $\gs = \grg$
there is nothing to prove. In the inductive proof we only need to look
at a place where the set increases. By the hypothesis on $\grg$ we can
suppose the result is true up to $\gs =\grg^-_{\ga_k, i}$ and try to
prove the result for $\gs =\grg^+_{\ga_k, i}$ (since new elements are
added only in these intervals). The new element added is $a^{\ga_k}_i$.
Denote this element by $a$. By hypothesis, $a/\GC_{\grg^-_{\ga_k, i}}$
\dns\ $N_{\ga_k}$ and so also not over $\bigcup_{k<n} N_{\ga_k} \cup
\bigcup_{k<n} (\GC_\grg \cap Y_{\ga_k})$. Now we can apply the induction
hypothesis and Proposition~\ref{trans}.

Notice that $X_\ga$ is contained in the definable closure of $Y_\ga$ and
vice versa so we also have $(\bigcup_{k < n} X_{\ga_k})/\GC_\grg$ \dns\
$\bigcup_{k<n} N_{\ga_k} \cup \bigcup_{k<n} (\GC_\grg \cap Y_{\ga_k})$.
We can immediately verify the induction hypothesis that if
$\grg<\delta$, $\grg\notin S$ then $\GC_\grg\slt \GC_{\gd +1}$. It is
enough to verify for $\ga_0, \dots, \ga_{n-1}$ and sufficiently large
$\gb$ that $(x_{\ga_0},\ldots, x_{\ga_{n-1}})/\GC_\gb$ does not split
over $\bigcup_{k<n} N_{\ga_k} \cup\bigcup_{k<n} (\GC_\gb \cap
Y_{\ga_k})$ (a set of size $<\gth$). But this sequence realizes the
ultrafilter average of $X_{\ga_0}\times \cdots \times X_{\ga_{n-1}}$. So
we are done by Proposition~\ref{ult}.

This completes the construction. Before continuing with the proof notice
that we get the following from the claim.

\begin{claim}
\label{nspl-claim}
1) For all $\ga \in \Ob$, $D \se \GC_{\gd+1}$, if $|D| < \gth$ then for
all but a bounded set of $i < \gth$, $D/\GC_{\grg_{\ga,i}^+}$ \dns\
$\GC_{\grg_{\ga,i}^-} \cup \{a^\ga_i\}$. Moreover for all but a bounded
set of $i < \gth$, $D/\GC_{\grg_{\ga,i}^+}$
\dns\ a subset of $\GC_{\grg_{\ga,i}^-} \cup \{a^\ga_i\}$ of size $<
\gth$.\\ 
2) for every subset $D$ of $\GC_{\gd+1}$ of cardinality $<\gth$ there is
a subset $w$ of $Ob(\gd$) of cardinality $<\gth$ and subset $Z$ of
$\GC_\gd$ of cardinality $<\gth$ such that: the type of $D$ over
$\GC_\gd$ does not split over $Z\cup
\bigcup_{j\in w} Y_j$.\\
3) In (2) for every large enough $i$, for every $\ga\in w$ the type of $D$
over $\GC_{\grg^+_{\ga, i}}$ does not split over $Z\cup 
\cup  \{Y_j\cap
\GC_{\grg^-_{\ga, i}}: j\in w\}\cup\{a_{\ga, i}\}$.\\
4) In (2), (3) we can allow $D\subseteq \GC$
\end{claim}

We now have to verify that $\GC$ has the desired properties.  Assume
that $f$ is an automorphism of $B=B^{\GC}$. We must show that

\begin{claim}
There is $\grg$ so that for all atoms $a$, $f(a)$ is definable with
parameters from $\GC_\grg$ and $a$.
\end{claim}

\proof (of the claim) Assume that $f$ is a counterexample. For every $\grg
\notin S$, choose an atom $a_\grg$ which witnesses the claim is false
with respect to $\GC_\grg$.  Since $\set{\gd < \gl}{\gd \notin S,
\cf \gd \geq \gth}$  is stationary, there is a set $N$ of cardinality
less than $\gth$ so that for a stationary set of $\grg$,
$a_\grg/\GC_\grg$ \dns\ $N$. In fact (since
$(\forall\ga<\gl)\ga^\gth<\gl$ as $\gl=\mu^+$, $\mu^{<\gth}=\mu)$ for
all but a bounded set of $\grg$ we can use the same $N$. Let $X$ code
the sequence $(\GC_\grg \colon \grg < \gl)$ and the function $f$.  Then,
by the previous discussion, $(\her{\rho}, \in, <, X)$ satisfies ``there
exists $N \se \GC$ so that $|N| < \gth$ and for all but a bounded set of
ordinals $\grg$, there is an atom $z$ so that $z/\GC_\grg$ \dns\
$N\mbox{''}$. Choose $\ga$ so that 
$$M^\ga_\gth\equiv_{M^\ga_\gth\cap\gl} (\her{\rho}, \in, <, X).$$

It is now easy to verify that an obstruction occurs at $\ga$.  Let $\gd
= \gz(\ga)$. In this case, $f_\ga$ is the restriction of $f$.  We use
the notation of the construction.  By the construction there is $D\se
\GC_{\gd+1}$, $\vert D\vert <\gth$ such that the type of $f(a_\ga$) over
$\GC_{\gd+1}$ does not split over $D$.  Apply 
Claim~\ref{nspl-claim}
above, parts (2),
(3) and get $Z
, w$ and $i^* < \gth$ (the $i^*$  is just explicating the 
``for every large enough $i$ to "for every $i \in [i^*,\gth)$'') .
  Let $D^*=Z\cup\cup\{Y_j \cap
C_{\grg^-_{\ga, i}}:j\in w\}$ so 
for every $i \in [i^*,\gth)$ we have

1) $f(a_\alpha)/\GC_{\gd+1}$ does not fork over $D (\subseteq
\GC_{\gd+1})$. \\
2) $D/\GC_{\gamma^+_{\alpha,i}}$ does not split over $D^*\cup a_{\alpha,i}$\\
3) $D^* \cup a_{\alpha,i}\subseteq \GC_{\gamma^+_{\alpha,i}}\subseteq
\GC_{\gd+1}$\\
so by the basic properties of non splitting
$f(x_\ga)/\GC_{\grg_{\ga, i}^+}$ \dns\ 
$D^*\cup \{a^\ga_i\}$, 
and note that we have : $D^* \cup N_\ga \se \GC_{\grg_{\ga, i}^-} $ and
$|D^* \cup N_\ga| < \gth$  where  $N_\ga$ comes from the construction .

An important point is that by elementariness for all ordinals $\tau \in
M^\ga_i \cap \gl$ and atoms $a \in M^\ga_i$ there is an ordinal $\gb$ so
that $\tau < \gb \in M^\ga_i \cap \gl$, $a \in B_\gb$, $\GC_\gb$ is
$\gth$-saturated (just take $\cf \gb \geq \gth$) and $f$ is an
automorphism of $B_\gb$. Choose such a $\gb \in M^\ga_{i+1}$ with
respect to $a^\ga_i$ and $\grg^-_{\ga,i}$
  
Since $f(a^\ga_i)$ is not definable from $a^\ga_i$ and parameters from
$\GC_{\grg_{\ga, i}^-}$ and $\GC_\gb$ is $\gth$-saturated there is $b
\neq f(a^\ga_i)$ realizing the same type over $D^* \cup \{a^\ga_i\} \cup
\set{f(a^\ga_j)}{j<i}$ with $b \in B_\gb$. Now for any atom $c \in
B_\gb$ we have (by the definition of $x_\ga$) that $c \leq x_\ga$ if
and only if $c = a^\ga_j$ for some $j \leq i$. Since this property is
preserved by $f$, we have that $f(a^\ga_i) \leq f(x_\ga)$ and $b \not
\leq f(x_\ga)$. But $\gb < \grg^+_{\ga, i}$ and
$f(x_\ga)/\GC_{\grg^+_{\ga, i}}$ \dns\ $\{a^\ga_i\} \cup D^*$.  So we
have arrived at a contradiction. \fin

In the proof above if we take $\gth$ to be uncountable then we can
strengthen the theorem (although we will not have any current use for
the stronger form).

\begin{theorem}
In the Theorem above if $\gth$ is uncountable then there is a finite set
of formulas $L'$ and a pseudofinite set $c$ so that for every atom $b
\in B$ $f(b)$ is $L'$-definable over $\{b\} \cup c$.
\end{theorem}

\proof  The argument so far has constructed a model in which
every automorphism of $B$ is pointwise definable on the atoms over some
$\GC_\grg$ (i.e., for every atom $b \in B$ $f(b)$ is definable from $b$
and parameters from $\GC_\grg$.  In the construction of the model we
have that every $\GC_\grg$ is contained in some pseudofinite set so we
are a long way towards our goal. To prove the theorem it remains to show
that we can restrict ourselves to a finite sublanguage. (Since all the
interpretations of the constants will be contained in $\GC_1$ we can
ignore them.)  Choose $\GC_\grg$ so that $f$ and $f^{-1}$ are pointwise
definable over $\GC_\grg$. Let $c$ be a pseudofinite set containing
$\GC_\grg$. We can assume that $f$ permutes the atoms of $B_\grg$.

Let the language $L$ be the union of an increasing chain of finite
sublanguages $(L_n \colon n < \go)$. Assume by way of contradiction that
for all $n$, $f$ is not pointwise definable on the atoms over any
pseudofinite set (and hence not over any $\GC_\ga$) using formulas from
$L_n$.  Choose a sequence $d_n$ of atoms so that for all $n$, $e_n =
f(d_n)$ is not $L_n$-definable over $\{d_n\} \cup c$ and both $d_{n+1}$
and $e_{n+1}$ are not definable over $\set{d_k}{k \leq n}\cup
\set{e_k}{k \leq n} \cup c$. Furthermore $d_{n+1}$ should not be
$L_{n}$-definable over $\set{d_k}{k \leq n} \cup c$.  The choice of
$d_n$ is possible by hypothesis, since only a pseudofinite set of
possibilities has been eliminated from the choice.

Let $\bar{d}, \bar{e}$ realize the average type (modulo some
ultrafilter) over $\GC_\grg \cup \{c\}$ of $\set{(d_k \colon k
\leq n), (e_k \colon k \leq n)}{n < \go}$.  These are pseudofinite
sequences which have $(d_n \colon n < \go)$ and $(e_n \colon n < \go)$
as initial segments. Say $\bar{d} = (d_i \colon i < n^*)$ for some
non-standard natural number $n^*$. Now let $x \in B$ be the join of
$\set{d_i}{i < n^*}$. (This join exists since $\bar{d}$ is a
pseudofinite sequence.)  For every $i < n^*$ there is a (standard)
$n_i$ so that $f(d_i)$ is $L_{n_i}$-definable over $\{d_i\} \cup
\GC_\grg$ and $d_i$ is $L_{n_i}$-definable over $\{f(d_i)\} \cup
\GC_\grg$. Since $\GC$ is $\gth$-saturated, the coinitiality
of $n^* \setminus \go$ is greater than $\go$. So there is some $n$ so
that $\set{i}{n_i = n}$ is coinitial. (Notice that for non-standard
$i$, there is no connection between $e_i$ and $f(d_i)$.)

Choose $k$ so that the formulas in $L_{n}$ have at most $k$ free
variables. Let $Z$ be the set of subsets $Y$ of $c$ of size $k$ so that
for all $ i < n^*$, neither $d_i$ nor $e_i$ is $L_{n}$-definable from $Y
\cup\set{d_j, e_j}{j < i}$. By the choice of $\bar{d}, \bar{e}$ every
subset of $\GC_\grg$ of size $k$ is an element of $Z$.
	
Consider $\go > m > n+1$. We claim there is no atom $y < f(x)$ such that
$d_m$ and $y$ are $L_n$ interdefinable over elements of $Z$.  Suppose
there is one and $y = f(d_i)$. Then $i$ is non-standard, since, by the
choice of $Z$, $y \neq e_j$ for all $j$. Since $d_i$ is definable from
$\{f(d_i)\} \cup c$, $d_i$ is definable from $\{d_m\} \cup c$.  This
contradicts the choice of the sequence. We can finally get our
contradiction. For $i < n^*$, we have $i < \go$ if and only if for all
$j > n +1$, if there is an atom $y < f(x)$ so that $d_j$ and $y$ are
$L_n$-interdefinable over elements of $Z$ then $i < j$. \fin

We will want to work a bit harder and get that the automorphisms in the
model above is actually definable on a large set.  To this end we prove
an easy graph theoretic lemma.

\begin{lemma}
\label{graph-lemma}
Suppose $G$ is a graph and there is $0 < k < \go$ so that the valence of
each vertex is at most $k$. Then there is a partition of (the set of
nodes of) $G$ into $k^2$ pieces $A_0, \ldots, A_{k^2-1}$ so that for any
$i$ and any node $v$, $v$ is adjacent to at most one element of $A_i$.
Furthermore if $\gl$ is an uncountable cardinal, each $A_i$ can be
chosen to meet any $\gl$ sets of cardinality $\gl$.  \end{lemma}

\proof Apply Zorn's lemma to get a sequence $\langle A_0, \ldots,
A_{k^2-1}\rangle$ of pairwise disjoint set of nodes such that for
$i<k^2$ and node $v$, $v$ is adjacent to at most one member of $A_i$ and
$\bigcup_{i< k^2} A_i$ maximal [under those constrains].  Suppose there
is $v$ which is not in any of the $A_i$. Since $v$ is not in any of the
$A_i$ for each $i$ there is $u_i$ adjacent to $v$ and $w_i\neq u_i$
adjacent to $v$ so that $w_i\in A_i$. If no such $u_i, w_i$ existed we
could extend the partition by adding $v$ to $A_i$. But as the valency of
every vertex of $G$ is $\leq k$ there are at most $k(k-1) < k^2$ such
pairs.

As for the second statement. Since the valence is finite every connected
component is countable. Hence we can partition the connected components
and then put them together to get a partition meeting every one of the
$\gl$ sets. \fin

Actually for infinite graphs we can get a sharper bound. Given $G$ we
form an associated graph by joining vertices if they have a common
neighbour. This gives a graph whose valence is at most $k^2 - k$. We
want to vertex colour this new graph. Obviously (see the proof above) it
can be vertex coloured in $k^2 +1 - k$ colours. In fact, by a theorem of
Brooks (\cite{Ore}, Theorem 6.5.1), the result can be sharpened further.
In our work we will only need that the colouring is finite, so these
sharpenings need not concern us.

We will want to form in $\GC$ a graph whose vertices are the atoms of
$B$ with a pseudofinite bound so that any atom $b$ is adjacent to
$f(b)$. This is easy to do in the case where all the definitions of
$f(b)$ from $\{b\} \cup c$ use only a finite sublanguage. In the general
case (i.e., when $\cf \gth$ may be $\go$) we have to cover the possible
definitions by a pseudofinite set.

\begin{lemma}
Continue with the notation of the proof. Then there is a pseudofinite
set $D$ and a pseudofinite natural number $k^*$ and a set of tuples $Z$
of length at most $k^*$ so that for every formula $\gf(\bar{x},\bar{y})$
there is $d \in D$ so that for all $\bar{a} \in c$ and $\bar{b} \in B$,
the tuple $(d, \bar{a}, \bar{b}) \in Z$ if and only if $\gf(\bar{a},
\bar{b})$.
\end{lemma}

\proof Since $\GC$ is $\ha_0$-saturated the lemma just says that a
certain type is consistent. Now $B \in \GC$ and $\GC$ satisfies the
separation scheme for all formulas (not just those of set theory).
Hence for any formula $\gf$, $\set{(\bar{a},\bar{b})}{\bar{a} \in c,
\bar{b} \in B\mbox{ and } \gf(\bar{a},\bar{b})}$ exists in $\GC$.  \fin 

Fix such sets $D$ and $Z$ for $c$. Say that an atom $a$ is {\em $D,
Z$-definable from} $b$ over $c$ if there are $d \in D$ and a tuple
(perhaps of non-standard length) $\bar{x} \in c$ so that $a$ is the
unique atom of $B$ so that $(d, \bar{x}, b, a) \in Z$. We say that $a$
and $b$ are {\em $D, Z$-interdefinable} over $c$ if $a$ is $D,
Z$-definable from $b$ over $c$ and $b$ is $D, Z$-definable from $a$ over
$c$. Notice (and this is the content of the last lemma) that if $a$ is
definable from $b$ over $c$ then it is $D, Z$-definable from $b$ over
$c$.

We continue now with the notation of the theorem and the model we have
built. Suppose $f$ is an automorphism which is pointwise definable over
a pseudofinite set $c$.  Define the {\em $c$-graph} to be the graph
whose vertices are the atoms of $B$ where $a, b$ are adjacent if $a$ and
$b$ are $D, Z$-interdefinable over $c$. Since $D$, $k^*$ and $c$ are
pseudofinite this is a graph whose valency is bounded by some
pseudofinite number. So in $\GC$ we will be able to apply
Lemma~\ref{graph-lemma}. We will say that $(A_i \colon i < n^*)$ is a
{\em good partition} of the $c$-graph if it is an element of $\GC$ which
partitions the atoms of $B$ into a pseudofinite number of pieces and for
any $i$ and any $a$, $a$ is adjacent to at most one element of $A_i$. If
$(A_i \colon i < n^*)$ is a good partition then for all $i, j$ let
$f_{i,j}$ be the partial function from $A_i$ to $A_j$ defined by letting
$f_{i,j}(a)$ be the unique element of $A_j$ if any so that $a$ is
adjacent to $f_{i,j}(a)$. Otherwise let $f_{i,j}(a)$ be undefined.

We have proved the following lemma.

\begin{lemma}
Use the notation above. For all $a \in A_i$ there is a unique $j$ so
that $f(a) = f_{i,j}(a)$
\end{lemma}

The proof of the last lemma applies in a more general setting.  Since we
will want to use it later, we will formulate a more general result.

\begin{lemma}
\label{cover-lemma}
Suppose $\GC$ is an $\ha_0$-saturated model of an expansion of
$\mbox{ZFC}^-$ which satisfies the separation scheme for all formulas.
If $A$, $B$ are sets in $\GC$ and $f$ is a bijection from $A$ to $B$ such
that $f$ and $f^{-1}$ are pointwise definable over some pseudofinite
set, then there is (in $\GC$) a partition of $A$ into pseudofinitely
many sets $(A_i \colon i < k^*)$ and a pseudofinite collection $(f_{i,j}
\colon i, j < k^*)$ of partial one-to-one functions so that for all $i$,
the domain of $f_{i,j}$ is contained in $A_i$ and for all $a \in A_i$
there exists $j$ so that $f_{i,j}(a) = f(a)$.  Moreover if we are given
(in $\GC$) a family $P$ of $|A|$ (in $\GC$'s sense) subsets of $A$ of
cardinality $|A|$, we can demand that for every $i<k^*$ and $A^*\in P$
we have $\vert A_i\cap A'\vert =\vert A\vert$ (in $\GC$'s sense).
\end{lemma}

\begin{lemma}
\label{definable}
Use the notation and assumptions above. For any $i<k^*$ we have: $f
\rest A_i$ is definable.
\end{lemma}

\proof  For each $ \grg < \gl$ such that $\grg \notin S$ and $\cf \grg
\geq \gth$ choose $y_\grg$ so that $y_\grg$ is the join of a
pseudofinite set of atoms contained in $A_i$ and containing $A_i \cap
B_\grg$.

\begin{claim}
For all $a \in A_i$ so that $a \leq y_\grg$, $f(a)$ is definable from
$f(y_\grg)$ and the set $\{f_{i,j}\colon j < n^*\}$. In particular this
claim applies to all $a \in A_i \cap B_\grg$.
\end{claim}

\proof (of the claim) First note that $f(a) < f(y_\grg)$ and so there is
$j$ such that $f_{i,j}(a) < f(y_\grg)$. Suppose that there is $k \neq j$
so that $f_{i, k}(a) < f(y_\grg)$. Choose $b \in A_i$ (not neccessarily
in $B_\grg$) so that $f(b) = f_{i, k}(a)$ (such a $b$ must exist since
every atom below $y_\grg$ is in $A_i$). But since $f_{i,k}$ is a partial
one-to-one function, $a = b$. But $f(a) = f_{i,j}(a) \neq f_{i,k}(a)$.
This is a contradiction, so $f \rest (A_i \cap B_\grg)$ is defined by
``$b = f(a)$ if there exists $g \in \set{f_{i,j}}{j < n^*}$ so that $b =
g(a)$ and $g(a) < f(y_\grg)$''.

Choose now $N_\grg$ of cardinality $< \gth$ so that the type of
$f(y_\grg), (f_{i,j} \colon j < n^*)/\GC_\grg$ \dns\ $N_\grg$. Notice
that $f \rest (A_i \cap B_\grg)$ is definable as a disjunction of types
over $N_\grg$, namely the types satisfied by the pairs $(a, f(a))$. By
Fodor's lemma and the cardinal arithmetic, there is $N$ and a stationary
set where all the $N_\grg = N$ and all the definitions as a disjunction
of types coincide. So we have that $f \rest A_i$ is defined as a
disjunction of types over $a$, $N$, set of cardinality $< \gth$. We now
want to improve the definability to definability by a formula. We show
that for all $\grg$, $f \rest (A_i \cap B_\grg)$ is defined by a formula
with parameters from $N$. This will suffice since some choice of
parameters and formula will work for unboundedly many (and hence all)
$\grg$.

Suppose that $f \rest (A_i \cap B_\grg)$ is not definable by a formula
with parameters from $N$. Consider the following type in variables $x_1,
x_2, z_1, z_2$.
\begin{quotation}
\noindent$\set{\gf(x_1, x_2) \mbox{ \iff\ }\gf(z_1, z_2)}{\gf \mbox{ a
formula with parameters from } N}$ $ \cup \{x_1, z_1 \leq y_\grg\}$$
\cup \{\psi(x_1, x_2), \neg \psi(z_1,z_2)\}$, 

\noindent where $\psi(u, v)$ is a formula saying `` there exists
$g \in \set{f_{i,j}}{j < n^*}$ so that $u = g(v)$ and $g(u) <
f(y_\grg)\mbox{''}$ 
\end{quotation}
This type is consistent since by hypothesis it is finitely satisfiable
in $A_i \cap B_\grg$. Since $\GC$ is $\gth$-saturated there are $a_1,
a_2, b_1, b_2$ realizing this type. Since $a_1, b_1 \leq y_\grg$
$\psi(a_1, a_2)$ implies that $f(a_1) = a_2$ and similarly $\neg
\psi(b_1, b_2)$ implies that $f(b_1) \neq b_2$. On the other hand
$a_1, a_2$ and $b_1, b_2$ realize the same type over $N$. This
contradicts the choice of $N$. \fin

The information we obtain above is a little less than we want since we
want a definablility requirement on every automorphism of every
definable Boolean algebra.  It is easy to modify the proof so that we
can get the stronger result.  Namely in the application of the Black Box
we attempt to predict the definition of a definable Boolean algebra in
the scope of $P$ as well as the construction of the model $\GC$ and an
automorphism. With this change the proof goes as before. So the
following stronger theorem is true.

\begin{theorem} 
There is a model $\GC$ of $T_1$ so that if $B \se P(\GC)$ is a definable
atomic Boolean algebra and $f$ is any automorphism of $B$ (as a Boolean
algebra) then there is a pseudofinite set $c$ such that for any atom $b
\in B$, $f(b)$ is definable from $\{b\} \cup c $.
\end{theorem}

As before if $\gth$ is taken to be uncountable then we can find a finite
sublanguage to use for all the definitions.

The proof of the analogue of the Theorem~\ref{Ba-thm} for ordered fields
is quite similar. Let $T_1$ denote the theory constructed in
Theorem~\ref{ofthy}.

\begin{theorem}
There is a model $\GC$ of $T_1$ so that if $F_1, F_2 \se P(\GC)$ are
definable ordered fields and $f$ is any isomorphism from $F_1$ to $F_2$
then there is a pseudofinite set $c$ so that for all $b \in F_1$ $f(b)$
is definable from $\{b\} \cup c$.
\end{theorem}

\proof To simplify the proof we will assume that $F = P(\GC)$ is an ordered
field and $f$ is an automorphism of $F$.  The construction is similar to
the one for Boolean algebras, although we will distinguish between the
case $\gth = \go$ and the case that $\gth$ is uncountable.  We will need
to take a little more care in the case of uncountable cofinalities.  The
difference from the case of Boolean algebras occurs at stages $\gd \in
S$. Again we predict the sequence of structures $(\GC_i\colon i < \gl)$
and an automorphism $f_\ga$ of the ordered field $F$. We say that {\em
an obstruction occurs at $\ga$ }if we can make the choices below. (The
intuition behind the definition of an obstruction is that there are
infinitesmals of arbitrarily high order for which the automorphism is
not definable.) There are two cases, the one where $\gth = \go$ and the
one where $\gth$ is uncountable.

Since it is somewhat simpler we will first consider the case $\gth =
\go$. Suppose an obstruction occurs at $\ga$. By this we mean that we
can choose $N_\ga \in M^\ga_0$ such that for each $i < \go$ we can
choose $a^\ga_i\in \gint{\ga}$ so that $a^\ga_i/\GC_{\grg^-_{\ga, i}}$
\dns\ $N_\ga$, $a^\ga_i >0$, $a^\ga_i$ makes the same cut in
$F_{\grg^-_{\ga, i}}$ as 0 does and $f(a^\ga_i)$ is not definable over
$a^\ga_i$ and parameters from $\GC_{\grg^-_{\ga, i}}$.

Then we let $x_i^\ga =\sum_{j \leq i} a^\ga_i$. We choose $x_\ga$ so
that it realizes the average type over some non-principal ultrafilter of
$\set{x_i^\ga}{i <\go}$. As before we can show that both the inductive
hypothesis and Claim~\ref{nspl-claim} are satisfied. Notice in the
construction that $f_\ga(x^\ga_i)$ is not definable over $\{a^\ga_i\}
\cup\GC_{\grg^-_{\ga,i}}$ , since $x^\ga_i = x^\ga_{i-1} + a^\ga_i$ and
$f(x^\ga_i) = f(x^\ga_{i-1}) + f(a^\ga_i)$ (here we conventionally let
$x^\ga_i = 0$ when it is undefined) and each $M^\ga_i$ is closed under
$f_\ga$. Notice as well that for all $i$, $x^\ga_i < x_{\ga}$ and
$x^\ga_i$ and $x_\ga$ make the same cut in $\GC_{\grg^+_{\ga,i}}$.

Suppose now that $f$ is an automorphism which is not pointwise definable
over any $\GC_\grg$. Since $F$ is an ordered field this is equivalent to
saying that $f$ is not definable on any interval.  Also since $\GC_\grg$
is contained in a pseudofinite set (in $\GC$) there is a positive
interval which makes the same cut in the pseudofinite set (and hence in
$F_\grg$) $0$ does. So we can choose $a_\grg$ in this interval so that
$f(a_\grg)$ is not definable over $\GC_\grg \cup \{a_\grg\}$.  Arguing
as before we can find a set $N \se \GC_{\grg}$ of cardinality less than
$\gth$ such that for all but boundedly many $\grg$ $a_\grg$ can be
chosen so that $a_\grg/\GC_\grg$ \dns\ $N$. Arguing as before we get an
ordinal $\ga$ so that an obstruction occurs at $\ga$.  Consider now
$f(x_\ga)$, $i$ and a finite set $N_\ga \se \GC_{\grg^-_{\ga, i}}$ so
that $f(x_\ga)/\GC_{\grg^+_{\ga,i}}$ \dns\ $N_\ga \cup \{a_i^\ga\}$.
Since $f(x^\ga_i)$ is not definable over $\{a^\ga_i\} \cup
\GC_{\grg_{\ga, i}}$, there is $b > f(x^\ga_i)$, $b \in F_{\grg^+_{\ga,
i}}$ which realizes the same type over $N_\ga \cup \{a^\ga_i\}$.
(Otherwise $f(x^\ga_i)$ would be the rightmost element satisfying some
formula with parameters in this set.) But $f(x^\ga_i) < f(x_\ga) <b$.
This contradicts the choice of $N_\ga$.

We now consider the case where $\gth$ is uncountable. Here we have to do
something different at limit ordinals $i$. If $i$ is a limit ordinal, we
choose $x^\ga_i\in \gint{\ga} \cap M^\ga_{i+1}$ to realize the average
type of the $\set{x^\ga_j}{j < i}$. For the successor case we let we let
$x^\ga_{i+1} = x^\ga_i + a_{i+1}^\ga$ (where the $a_{j}^\ga$ are chosen
as before). In this construction for all $\grg$ and all $i$
$x_i^\ga/\GC_\grg$ \dns\ $N_\ga \cup \set{x^\ga_j}{\grg^-_{\ga, j} <
\grg }$. This is enough to verify the inductive hypothesis and
Claim~\ref{nspl-claim} if we restrict the statement to a successor
ordinal $i$. The rest of the proof can be finished as above.\fin

To apply the theorem we explicate the explanation in the introduction.
\begin{Observation}
\label{4.x}
For proving the compactness of $L(\Qof)$ (or $L(Q_{BA})$ etc) it
suffices to prove
\item[($*$)] Let $T$ be a first order theory $T$ and let $P$, $R$ be
unary predicates such that for every model $M$ of $T$ and an
automorphism $f$ of a definable ordinal field (or Boolean algebra etc)
$\se P^M$ definable (with parameter), for some $c\in M$, $R(x, y, c)$
defines $f$.

{\em Then} $T$ has a model $M^*$ such that: every automorphism $f$ of
definable ordered field (or Boolean algebra etc.) $\se P^M$, $f$ is
defined in $M^*$ by $R(x, y, c)$ for some $c\in M$.
\end{Observation}

\proof Assume $(*)$ and let $T_0$ be a given theory in the stronger logic;
without loss of generality all formulas are equivalent to relations.
For simplicity ignore function symbols.

For every model $M$ let $M'$ be the model with the universe $\vert M\vert
\cup \{f:f$ a partial function from $\vert M\vert$ to
$\vert M\vert\}$, and relations those of $M$, $P^{M'}=\vert M\vert$,
$R^{M'}=\{(f, a, b):a, b\in M$, $f$ a partial function from $\vert
M\vert$ to $\vert M\vert$ and $f(a)=b\}$.  There is a parallel
definition of $T'$ from $T$, and if we intersect with the first order
logic, we get a theory to which it suffices to apply $(*)$.

\section{Augmented Boolean Algebras}
\setcounter{theorem}{0}

For Boolean algebras we don't have the full result we would like to have
but we can define the notion of augmented Boolean algebras and then
prove the compactness theorem for quantification over automorphisms of
these structures.

\medskip

\noindent {\sc Definition.} An {\em augmented  Boolean algebra} is a
structure $(B, \leq, I, P)$ so that $(B, \leq)$ is an atomic Boolean
algebra, $I$ is an ideal of $B$ which contains all the atoms, $P \se B$,
$|P| > 1$, for $x \neq y \in P$ the symmetric difference of $x$ and $y$
is not in $I$, and for all atoms $x \neq y$ there is $z\in P$ such that
either $x \leq z$ and $y \not\leq z$ or vice versa.

\medskip 

Notice that if we know the restriction of an automorphism $f$ of an
augmented Boolean algebra to $P$ then we can recover $f$. Since for any
atom $x$, $f(x)$ is the unique $z$ so that for all $y\in P$ $z
\leq f(y)$ if and only if $x \leq y$ and the action of $f$ on
the atoms determine its action on the whole Boolean algebra.

Let $Q_{\rm Aug}$ be the quantifier whose interpretation is ``$Q_{\rm
Aug} f\,(B_1, B_2) \ldots$'' holds if there is an isomorphism $f$ of the
augmented Boolean algebras $B_1,B_2$ so that~$\ldots$.

\begin{theorem}
The logic {\rm  L($Q_{\rm Aug}$)} is compact.
\end{theorem}

\proof By 3.1 (and 4.9.1) it is enough to show that if $T$ is a theory
which says that every automorphism of a definable Boolean algebra
$\subseteq P$ is definable by a fixed formula and $T_1$ and $\GC$ are as
above then every automorphism $f$ of a definable augmented Boolean
algebra $(B, \leq, I, P)$ of $P(\GC)$ is definable. We work for the
moment in $\GC$. First note that $I$ contains the pseudofinite sets so
$\GC$ thinks that the cardinality of $B$ is some $\mu^*$ and for every
$x, y \in P$ the symmetric difference of $x$ and $y$ contains $\mu^*$
atoms (since $\GC$ thinks that $B$ is $\mu^*$ saturated).

So, by Lemma~\ref{cover-lemma} and Lemma~\ref{definable}, we can find a
set of atoms $A$ so that $f \rest A$ is definable and for every $x$ and
$y$, $A$ contains some element in the symmetric difference. Now we can
define $f(y) $ as the unique $z \in B$ so that for all $a \in A$, $a
\leq y$ if and only if $f(a) \leq z$. Clearly $f(y)$ has this property.
Suppose for the moment that there is some $z \neq f(y)$ which also has
this property.  Since $f$ is an automorphism there is $x \in P$ so that
$f(x) = z$.  Now choose $a\in A$ in the symmetric difference of $x$ and
$y$. Then $f(a)$ is in exactly one of $f(y)$ and $f(x)$. \fin

\section{Ordered Fields}
\setcounter{theorem}{0}

Here we will prove the compactness of the quantifier $ Q_{\rm Of}$. Let
$Q_{\rm Of}$ be the quantifier whose interpretation is ``$Q_{\rm Of}
f\,(F_1, F_2) \ldots$'' holds if there is an isomorphism $f$ of the
ordered fields $F_1,F_2$ so that $\ldots$.  We will use various facts
about dense linear orders. A subset of a dense linear order is {\em
\swd\ }if it is dense in some non-empty open interval. A subset which is
not \swd\ is {\em \nwd}. The first few properties are standard and
follow easily from the fact that a finite union of \nwd\ sets is \nwd.

\begin{proposition}
If a \swd\ set is divided into finitely many pieces one of the pieces is
\swd.
\end{proposition}


\begin{proposition}
If $\set{f_k}{k \in K}$ is a finite set of partial functions defined on
a \swd\ set $A$ then there is a \swd\ subset $A' \se A$ so that each
$f_k$ is either total on $A'$ or no element of the domain of $f_k$ is in
$A'$.
\end{proposition}

We can define an equivalence relation on functions to a linearly ordered
set by $f\equiv g$ if for every non-empty open interval $I$ the
symmetric difference of $f^{-1}(I)$ and $g^{-1}(I)$ is \nwd.

\begin{proposition}
\label{part-prop}
Suppose $\cal F$ is a finite collection of functions from a \swd\ set
$A$ to a linearly ordered set. Then there is $\cal I$ a collection of
disjoint intervals (not necessarily open) and a \swd\ set $A' \se A$ so
that for all $f, g \in \cal F$ either $f \rest A' \equiv g \rest A'$ or
there are disjoint intervals $I, J\in \cal I$ so that $f(A') \se I$ and
$g(A') \se J$.
\end{proposition}

\proof The proof is by induction on the cardinality of $\cal F$.
Suppose that ${\cal F} = {\cal G}\cup\{g\}$ and $\cal J$ is a set of
intervals and $A'' \se A$ is a \swd\ set satisfying the conclusion of
the theorem with respect to $\cal G$. If there are some $f \in \cal G$
and $A'''$ a \swd\ subset of $A''$ such that $f \rest A''' \equiv g
\rest A'''$ then we can choose a \swd\ $A'\subseteq A'''$ such that for
all $I\in {\cal J}$ we have $(f\rest A')^{-1}(I)=(g\rest A')^{-1}(I)$.
Then $A'$ and $\cal J$ are as required in the theorem. So we can assume
that no such $f$ and $A'''$ exist.

Let $\{f_0, \ldots, f_n\}$ enumerate a set of $\equiv$ class
representatives of $\cal G$ and let ${\cal J} = \{J_0, \ldots, J_n\}$
where $f_k(A'') \se J_k$ for all $k$ (so for $\ell<m\le n$ $J_ell\cap
J_m=\emptyset)$. By induction on $k$ we will define a descending
sequence $A_k$ of \swd\ subset of $A''$ and intervals $I_k \se J_k$ with
the property that $f_k(A_k) \se I_k$ and $g(A_k)$ is disjoint from
$I_k$.  Let $A_{-1} = A''$.  Consider any $k$ and suppose $A_{k-1}$ has
been defined. If possible, choose an interval $J_k' \se J_k$ so that the
symmetric difference of $(g\rest A_{k-1})^{-1}(J_k')$ and $(f_k\rest
A_{k-1})^{-1}(J_k')$ is \swd.  There are two possibilities. If
$(f_k\rest A_{k-1})^{-1}(J_k') \setminus (g\rest A_{k-1})^{-1}(J_k')$ is
\swd, then let $A_k = (f_k\rest A_{k-1})^{-1}(J_k') \setminus (g\rest
A_{k-1})^{-1}(J_k')$ and let $I_k = J_k'$. Otherwise we can choose $I_k$
a subinterval of $J_k \setminus J_k'$ and $A_k$ a \swd\ subset of
$(g\rest A_{k-1})^{-1}(J_k')$ so that $f_k(A_k) \se I_k$.  If there was
no respective $J_k'$ then we would have $f_k\rest A_{k-1}\equiv g\rest
A_{k-1}$ - contrary to our assumption. Given $A_n$ and $I_0, \ldots,
I_n$, we can choose $A_{n+1} \se A_n$ \swd\ and an interval $I_{n+1}$
disjoint from $I_k$, for all $k$ so that $g(A_{n+1}) \se I_{n+1}$.
Finally we let $A' = A_{n+1}\cap \bigcap_{f\in \cal G}
f^{-1}(\bigcup_{k\leq n} I_k)$ and ${\cal I}=\{I_0,\ldots,I_n\}$. 
\fin

\begin{theorem}
The logic ${\rm L}(Q_{\rm Of})$ is compact.
\end{theorem}

\proof We work in the model $\GC$ constructed before (suffice
by 3.2 (and 4.9.1)). Suppose we have two definable ordered fields $F_1,
F_2$ contained in $P(\GC)$ and an isomorphism $f$ between them. Since
$f$ and $f^{-1}$ are locally definable over a pseudofinite set there are
$(A_i \colon i < k^*)$ and $(f_{i, j} \colon i, j < k^*)$ as in
Lemma~\ref{cover-lemma}. Since $\GC$ satisfies that $k^*$ is finite,
there is some $i$ so that $A_i$ is \swd. Fix such an $i$ and for
notational simplicity drop the subscript $i$. So $f_j$ denotes
$f_{i,j}$. By restricting to a \swd\ subset (and perhaps eliminating
some of the $f_j$) we can assume that each $f_j$ is a total function (we
work in $\GC$ to make this choice) (use 6.2). Again choosing a \swd\
subset we can find a \swd\ set $A$ and a collection $\cal I$ of disjoint
intervals of $F_2$ so that the conclusion of Proposition~\ref{part-prop}
is satisfied. Since each $f_j$ is one-to-one we can assume that all the
intervals in $\cal I$ are open. Choose now an interval $J$ contained in
$F_1$ so that $A$ is dense in $J$.

Next choose $a_1 < a_2 \in J$ and an interval $I \in \cal I$ so that
$f(a_1), f(a_2) \in I$. To see that such objects exist, first choose any
$a_1 \in A \cap J$. There is a unique interval $I \in \cal I$ so that
$f(a_1) \in I$.  Choose $b \in I$ so that $f(a_1) < b $ and let $a_2$ be
any element of $J \cap A$ such that $a_1 < a_2 < f^{-1}(b)$. Of course
this definition of $a_1, a_2$ cannot be made in $\GC$, but the triple
$(a_1, a_2, I)$ exists in $\GC$. Without loss of generality we can
assume that $A \se (a_1, a_2) = J$. Consider now ${\cal F} =
\set{f_j}{f_j(A) \se I}$. This is an equivalence class of functions and
for every $a \in A$, there is some $g \in \cal F$ so that $g(a) = f(a)$.

The crucial fact is that for all $b \in J$ and $g \in {\cal F}$, $B =
\set{a \in A}{ b < a \mbox{ and } g(a) < f(b)}$ is \nwd.  Assume not.
Since $\cal F$ is an equivalence class, for all $h \in
\cal F$ $\set{a \in B}{h(a) \geq f(b)}$ is a \nwd\  set. But since
$\cal F$ is a pseudofinite set, $\set{a \in B}{\mbox{for some } h \in
{\cal F}, h(a)\geq f(b)}$ is \nwd. So there is some $a \in A$ so that $b
< a$ and for all $g \in \cal F$ $g(a) < f(b)$. This gives a
contradiction since there is some $g$ so that $f(a) = g(a)$. Similarly,
$\set{a \in A}{a < b\mbox{ and }g(a) > f(b)}$ is \nwd.

	With this fact in hand we can define $f$ on $J$, by $f(b)$ is
the greatest $x$ so that for all $g \in \cal F$, both $\set{a \in
A}{b < a \mbox{ and }g(a) < x}$ and $\set{a \in A}{a < b \mbox{ and
}g(a) > x}$ are \nwd. Since an 
isomorphism between ordered fields is definable if and only if it is
definable on an interval we are done. \fin

\end{document}